\pdfoutput=1 
\documentclass[11pt,oneside]{amsart}
\usepackage[leqno]{amsmath}
\usepackage{amssymb,amsthm} 
\usepackage[parfill]{parskip} 
\usepackage{libertine} 
\usepackage[T1]{fontenc}
\usepackage{textcomp}
\usepackage[varqu,varl]{inconsolata} 
\usepackage{amsmath,amsthm}
\usepackage[libertine,bigdelims,vvarbb]{newtxmath}
\usepackage{graphicx} 
\usepackage{calc}

\newcommand{\beq}{\begin{equation}}

\newcommand{\eeq}{\end{equation}}
\numberwithin{equation}{section} 
\newtheorem{thm}[equation]{Theorem}

\newtheorem{cor}[equation]{Corollary}

\theoremstyle{remark}

\theoremstyle{definition}
\newtheorem{defn}[equation]{Definition}

\title{On P=NP  either false or independent of ZFC}

\usepackage[paperwidth=7in, paperheight=10in, textwidth=4.75in, textheight=8in, includehead=true]{geometry}
\usepackage{verbatim} 
\usepackage{url}
\usepackage[pdfborderstyle={/S/U/W 1}]{hyperref}

\def\image{\operatorname{image}}

\let\emptyset\emptysetAlt
\hypersetup{
pdftitle={ZFC},
pdfauthor={\textcopyright\ S. G. Williamson}
}
\author{S. Gill Williamson}

\thanks{Computer Science and Engineering, 
University of California San Diego; \url{http://cseweb.ucsd.edu/~gill/} or search
"s gill williamson home" {\bf Keywords:}  ZFC independence, P equals NP, ZFC Limbo
 regressive regularity, subset sum problem, Jump Free Theorem}

\begin{document}
\begin{abstract}
  
Our main result, Theorem~\ref{thm:polyprop}, uses  Friedman's Jump Free Theorem, Theorem~\ref{thm:jumpfree}, 
which he has shown to be independent of  ZFC, the usual axioms of set theory~\cite{hf:alc}. 
 We conjecture that Theorem~\ref{thm:polyprop}, a straight forward translation of the statement of Theorem~\ref{thm:jumpfree} into sets and functions, is also independent of ZFC as is its immediate Corollary~\ref{cor:subsum}.
 It is easy to show that a proof that $P=NP$  will prove Corollary~\ref{cor:subsum}.  
If Corollary~\ref{cor:subsum} is in fact independent of ZFC then a  ZFC proof of $P=NP$ is impossible and, in fact, $P=NP$ is  itself independent of ZFC (or perhaps false).  

\end{abstract}

 \maketitle
 

\section{Introduction}  Reference~\cite{hf:nlc} presented for the first time numerous examples of statements in discrete and finite mathematics which are proved using large cardinals and shown to require them.  Reference~\cite{hf:alc} evolved from a seminar presented at UCSD (University of California San Diego) by Friedman and subsequent discussions.  This latter work has inspired the studies of combinatorial structures cited in the references at the end of this paper. The paper in the Journal of Combinatorics ~\cite{gw:lim} describes connections between combinatorics, ZFC independence and the subset sum problem. An  exposition of how this connection is made is given in ~\cite{gw:adn}. An alternative approach is presented in this short note.


\section{Basic definitions and theorems}

$Z$ denotes the integers and  $N$ the nonnegative integers. 
For $x=(n_1, \ldots, n_k)\in N^k$, $\max\{n_i\mid i=1,\ldots, k\}$ is
denoted by $\max(x)$.  Define $\min(x)$ similarly. 

\begin{defn}[\bf Cubes and Cartesian powers in $N^k$]
The set  $E_1\times\cdots\times E_k$, where $E_i\subset N$, $|E_i|=p$, $i=1,\ldots, k,$   is called a $k$-cube of length $p$.  If $E_i = E, i=1,\ldots, k,$ then this cube is  $E^k$, the $k$th Cartesian power of $E$.
\label{def:cubespowers}
\end{defn}

\begin{defn}[\bfseries Order equivalent $k$-tuples]
\label{def:ordtypeqv}
Two $k$-tuples, $x=(n_1,\ldots,n_k)$ and $y=(m_1,\ldots,m_k)$, are  
{\em order equivalent tuples $(x\, ot \,y)$} if the following holds:
$\{(i,j)\mid n_i < n_j\} =  \{(i,j)\mid m_i < m_j\}$ and  $\{(i,j)\mid n_i = n_j\} =  \{(i,j)\mid m_i =m_j\}.$  
\end{defn}

Note that $ot$ is  an equivalence relation on $N^k$.
We use ``$x\,ot\,y$'' and ``$x,\,y$ of order type $ot$'' to mean $x$ and $y$ belong to the same order type equivalence class.  The number of equivalence classes is bounded by $k^k$.

We present some basic definitions due to Friedman~\cite{hf:alc}, \cite{hf:nlc}.


\begin{defn}[\bf Field of a function and reflexive functions]
\label{def:fldreflfncs}
For $A\subseteq N^k$ define ${\rm field}(A)$ to be the set of all coordinates of elements of $A$.  A function $f$ is {\em reflexive} if 
${\rm domain}(f) \subseteq N^k$ and  ${\rm range}(f) \subseteq {\rm field}({\rm domain}(f))$.
\end{defn}

\begin{defn}[{\bf The set of functions} $T(k)$ ]
\label{def:tk}
$T(k)$ denotes all reflexive functions with finite domain.  We denote a function with domain $D\subseteq N^k$ by $f_D$.
\end{defn}

\begin{defn} [\bf Full and jump free families]
\label{def:fulljump} 
Let $Q\subseteq T(k).$
\begin{enumerate}

\item {\bf full family:}  We say that $Q\subseteq T(k)$ is a {\em full} family of functions if for every finite subset 
$D\subset N^k$ there is at least one function $f$ in $Q$ with domain  $D$.

\item{\bf jump free family:} For any finite $D\subset N^k$ and for any $x\in D$ we define $D_x = \{z\mid z\in D,\, \max(z) < \max(x)\}$. 
Suppose that for all functions $f_A$ and $f_B$  in $Q$,  
 $x\in A\cap B$, $A_x \subseteq B_x$, and $f_A(y) = f_B(y)$ for all $y\in A_x$ imply that 
$f_A(x) \geq f_B(x)$.  Then $Q$ will be called a {\em jump free} family of functions in $N^k$.  
\end{enumerate}
\label{def:fullrefjf}
\end{defn} 

\begin{figure}[h]
\begin{center}
\includegraphics[scale=.85]{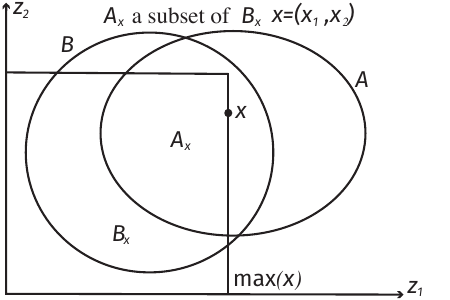}
\caption{Basic jump free condition~\ref{def:fullrefjf}.}
\label{fig:jfvenn}
\end{center}
\end{figure}

The jump free property arises from numerous recursively constructed algorithms in combinatorics. 
See \cite{hf:alc}, \cite{gw:lem}, \cite{gw:sub}, \cite{gw:lim}. 
We use ZFC for the axioms of set theory: Zermelo-Frankel  plus the axiom of choice. 
Freedman's  Jump Free Theorem, Theorem~\ref{thm:jumpfree} stated below, can be proved in 
 ZFC + ($\forall n$)($\exists$ $n$-subtle cardinal)
 but not in ZFC+ ($\exists$ $n$-subtle cardinal) for any fixed $n$ (assuming this theory is consistent).
\begin{defn}[\bf Function regressively regular over $E^k$]
\label{def:regreg}
Let $k\geq 2$, $D\subset N^k$, $D$ finite, $f: D\rightarrow N$. 
We say that $f$ is {\em regressively regular} over $E^k$, 
$E^k\subseteq D$, if for each order type equivalence class $ot$ of $k$-tuples of $E^k$ either (1) or (2) occurs:  (1) For all $x, y\,\in E^k$ of order type $ot$, $f(x)=f(y)< \min(E)$ or (2) For all $x\in E^k$ of order type $ot$, $f(x)\geq \min(x)$.  
\end{defn}


\begin{thm}[\bf Jump Free Theorem (\cite{hf:alc}, \cite{hf:nlc})] 
\label{thm:jumpfree}
Let $S\subseteq T(k)$ be  a full and jump free family of functions. Let $p, k\geq 2$. 
Then some $f\in S$ is regressively regular over some $E^k$, cardinality $|E|=p$.
\end{thm}
A proof of the Jump Free Theorem is in Section 2 of \cite{hf:alc},
``Applications of Large Cardinals to Graph Theory.''  
A discussion of a vector valued extension of the Jump Free Theorem and some of its applications is given in Section 3 of \cite{jg:pos} as well as a short introduction to n-subtle cardinals (Appendix A).

\section{Subset sum instances and the Jump Free Theorem}

\begin{defn}
\label{def:fncblocks}({\bf Partitions of N})
Let $S \subseteq T(k)$ be  a full and jump free family of functions.
For a function  $f \in S$ and $x \in E^k\subseteq domain(f)$, $f(x)$ is in one of three intervals
that form a partition of $N$:
 $I^k_0=[0, \min(E))$,  $I^k_1=[\min(E), \min(x))$,  $I^k_2=[\min(x),\infty)$. 
\end{defn}

Note that for $x\in E^k$,
 $f(x) \notin [\min(E), \min(x))$ if $f$ is {\em regressively regular} over $E^k$. 
 In general, we define sets of integers in $Z$ as follows. 
 
 Let $\Gamma=\{\gamma\mid \gamma: N  \rightarrow Z\}$ denote bijections from $N$ to $Z$.  Let $\gamma X = \{\gamma(x)\mid x\in X\}$, $\gamma\emptyset = \emptyset$. Let $\hat f = f | E^k$ be $f$ restricted to $E^k$.
 Note that 
 $$\image(\hat f)\cap I^k_i = \{\hat f(x)\mid x \in I^k_i\}, i=0, 1, 2.$$
 We prefer the former notation for what follows.
 
\begin{defn}
\label{def:collections of sets}
 Define  subsets of $Z$ as follows :  
For  $\gamma_0$, $\gamma_1$, $\gamma_2$ in $\Gamma$, $f\in S$, $E^k\subseteq domain(f)$, $|E|=p$,  $\hat f = f | E^k$.  Define
$$F(\hat f)= \gamma_0(\image(\hat f)\cap I^k_0) \cup \gamma_1 ( \image(\hat f)\cap I^k_1) \cup  \gamma_2(\image(\hat f)\cap I^k_2)$$
 and 
$$H(\hat f) =\gamma_0(\image(\hat f)\cap I^k_0) \cup \gamma_2(\image(\hat f)\cap I^k_2).$$

\end{defn}

\begin{thm}{\bf (Jump Free Theorem Set Version).}
\label{thm:polyprop}
Let $S\subseteq T(k)$, $k\geq 2$, be  a full and jump free family of functions. 
Consider  sets $F(\hat f)$  and $H(\hat f)$ of integers of Definition~\ref{def:collections of sets}.  For each $p\geq 2$ there exists $f_p\in S\subseteq T(k)$  for which  $F(\hat f_p) = H(\hat f_p)$.
\begin{proof}
For each $p \geq 2$ use the Jump Free Theorem~\ref{thm:jumpfree} to choose $ f_p\in S$  regressively regular over some $E^k \subseteq {\mathrm domain(f_p)}$, $|E| = p$. By regressive regularity $\gamma_1 (\image(\hat f_p)\cap I^k_1)$ is the empty set. 
Thus, $F(\hat f_p)= H(\hat f_p)$ for this $f_p$. 
\end{proof} 
\end{thm}
We conjecture that  Theorem~\ref{thm:polyprop} like the Jump Free Theorem is independent of ZFC.     Theorem~\ref{thm:polyprop} has been proved by the ZFC independent Jump Free Theorem~\ref{thm:jumpfree} and no other proof is evident.
A corollary of Theorem~\ref{thm:polyprop} is as follows

\begin{cor}{\bf( Subset Sum)}
\label{cor:subsum}
There exists $t \in N$ such that for
 $p\geq 2$ there exists $f_p\in S\subseteq T(k)$  for which  $F(\hat f_p)$ is subset sum target zero solvable in polynomial time $O(p^{kt})$ if and only if 
$H(\hat f_p)$ is subset sum target zero is solvable in polynomial time $O(p^{kt})$. Here the cardinalities of $F(\hat f_p)$ and $H(\hat f_p)$  are $p^k$ as multisets. 
\begin{proof}
Follows directly from Theorem~\ref{thm:polyprop}.  It also follows if $P=NP$.
\end{proof}
\end{cor}
It is possible that this corollary itself might be independent of ZFC. 
In which case this would imply that $P=NP$ itself is independent of ZFC or perhaps false. 
~\cite{gw:adn}

%

\bibliographystyle{alpha}
\bibliography{jeffjoc}

\end{document}